\documentclass[12pt]{article}

\advance \topmargin by -\headheight
\advance \topmargin by -\headsep
     
\evensidemargin \oddsidemargin
\usepackage{amssymb,latexsym,pstricks}

\newtheorem{thm}{Theorem}[section]
\newtheorem{lem}[thm]{Lemma}
\newtheorem{cor}[thm]{Corollary}

\newcommand{\bbox}{\hspace{\stretch{1}} $\blacksquare$}

\begin{document}

\bibliographystyle{elsart-num-sort}

\title{Permutation Statistics and $q$-Fibonacci Numbers}

\author{Adam M. Goyt\footnote{email: goytadam@mnstate.edu} and David Mathisen\footnote{email: mathisda@mnstate.edu}\\ Mathematics Department\\ Minnesota State University Moorhead\\ Moorhead, MN 56562}

\maketitle
\begin{abstract}
In a recent paper, Goyt and Sagan studied distributions of certain set partition statistics over pattern restricted sets of set partitions that were counted by the Fibonacci numbers.  Their study produced a class of $q$-Fibonacci numbers, which they related to $q$-Fibonacci numbers studied by Carlitz and Cigler.  In this paper we will study the distributions of some Mahonian statistics over pattern restricted sets of permutations.  We will give bijective proofs connecting some of our $q$-Fibonacci numbers to those of Carlitz, Cigler, Goyt and Sagan.  We encode these permutations as words and use a weight to produce bijective proofs of $q$-Fibonacci identities.  Finally, we study the distribution of some of these statistics on pattern restricted permutations that West showed were counted by even Fibonacci numbers.
\end{abstract} 

\section{Introduction}

We will study the distribution of two Mahonian statistics, {\it inv} and {\it maj}, over sets of pattern restricted permutations.  In particular, we will study the distributions of these statistics over pattern-restricted sets which are counted by the Fibonacci numbers.  These distributions will give us $q$-Fibonacci numbers which are related to the $q$-Fibonacci numbers of Carlitz~\cite{Carlitzfibnotes3}, Cigler~\cite{Ciglerfib2,Ciglerqfib1}, and Goyt and Sagan~\cite{SagGoytqFib}. 

Let the $n^{th}$ Fibonacci number be $F_n$, where $F_n=F_{n-1}+F_{n-2}$ and $F_0=1$ and $F_1=1$.  Let $[n]=\{1,2,\dots n\}$ and $[k,n]=\{k,k+1,\dots, n\}$.  We will call two integer sequences $a_1a_2\dots a_k$ and $b_1b_2\dots b_k$ are {\it order isomorphic} if $a_i<a_j$ whenever $b_i<b_j$.  Let $S_n$ be the set of permutations of $[n]$, and suppose $\pi=p_1p_2\dots p_m\in S_m$ and $\sigma=q_1q_2\dots q_n\in S_n$.  We say that $\sigma$ contains the {\it pattern} $\pi$ if there is a subsequence $\sigma'=q_{i_1}q_{i_2}\dots q_{i_m}$ of $\sigma$ which is order isomorphic to $\pi$, otherwise we say that $\sigma$ {\it avoids} $\pi$.  For example, a copy of $\pi=321$ in $\sigma=564312$ is $641$.  However, $\sigma$ avoids $123$ because it does not have any increasing subsequences of length three.  Let $R$ be a set of patterns and let $S_n(R)$ be the set of permutations in $S_n$ that avoid every pattern in $R$.

The sets that we wish to study are $S_n(123,132,213)$ and $S_n(231,312,321)$.  To study statistical distributions on these sets it will be necessary to understand their structure.  For two sequences of integers $\alpha$ and $\beta$, we will say $\alpha<\beta$ if $\max \alpha<\min \beta$.  A permutation is called {\it layered} if it can be written $\pi=\pi_1\pi_2\dots\pi_k$ where $\pi_i<\pi_j$ whenever $i<j$, and the $\pi_i$ are decreasing.  The $\pi_i$ will be called {\it layers}.  For example, $\pi=321549876$ is layered with layers $321$, $54$, and $9876$.       

For a permutation $\pi=p_1p_2\dots p_n$, let the {\it reversal} of $\pi$ be $\bar{\pi}=p_np_{n-1}\dots p_1$.  The reversal of $\pi$ given above is $\bar{\pi}=678945123$.  If $\pi$ is layered, the $\bar{\pi}$ is reverse layered and $\bar{\pi_i}$ is a layer of $\bar{\pi}$ whenever $\pi_i$ is a layer of $\pi$.  We will call a (reverse) layered permutation a {\it matching} if all layers of the permutation are of size at most two.  For example, $\pi=6753421$ is a reverse layered matching.  A layer with one element is a {\it singleton}, and a layer with two elements is a {\it doubleton}.

\begin{thm}  $S_n(123,132,213)$ is the set of reverse layered matchings of $[n]$.\end{thm}

{\it Proof:}  It is clear from the definition of a reverse layered matching that reverse layered matchings avoid $123$, $132$ and $213$.  

Since $S_0(123,132,213)$ contains only the empty permutation and $S_1(123,132,213)=\{1\}$, these two sets consist entirely of reverse layered matchings.  

Let $\pi=p_1p_2\dots p_n\in S_n(123,132,213)$.  If $p_1\not=n$ and $p_2\not=n$ then there are two elements to the left of $n$ in $\pi$, so there is either a copy of 123 or 213.  Thus, $p_1=n$ or $p_2=n$.  

If $p_1=n$, then $\pi=np_2\dots p_n$.  If $p_2=n$ and $p_1\not=n-1$, then $\pi$ contains a copy of $132$.  Thus $\pi=(n-1)np_3\dots p_n$.  In either case $\pi$ is a reverse layered matching.  \bbox

\medskip

The following is an immediate consequence of Theorem 1.1.

\begin{cor}  $S_n(231,312,321)$ is the set of layered matchings of $[n]$.\end{cor} \bbox

\medskip

We will focus on the distributions of the two Mahonian statistics, {\it maj} and {\it inv}, over $S_n(123,132,213)$ and $S_n(231,312,321)$.  If $\pi=p_1p_2\dots p_n$, then an {\it inversion} is any pair $p_i,p_j$ where $i<j$ and $p_i>p_j$.  We define $inv(\pi)$ to be the number of inversions in $\pi$.  A {\it descent} in $\pi$ is a pair $p_ip_{i+1}$ such that $p_i>p_{i+1}$.  Let $D(\pi)=\{i:p_ip_{i+1}\mbox{ is a descent}\}$, then we define the major index to be $$maj(\pi)=\sum_{i\in D(\pi)}i.$$

As we mentioned before, we are interested in studying the distributions of these statistics over the sets $S_n(123,132,213)$ and $S_n(231,312,321)$.  It was shown by Simion and Schmidt~\cite{SimionSchmidt} that $|S_n(123,132,213)|=|S_n(231,312,321)|=F_n$.  Thus, each distribution will give us a $q$-analogue of the Fibonacci numbers ($q$-Fibonacci numbers).  

In the next section we will encode these permutations as words and define a weight that will give the $q$-Fibonacci numbers that we are interested in.  In Section 3, we will give bijective proofs that the two $q$-Fibonacci numbers produced by the distribution of $maj$ are the same as those studied by Cigler~\cite{Ciglerqfib1} and Goyt and Sagan~\cite{SagGoytqFib}.  In Section 4, we will use the techniques developed by Benjamin and Quinn~\cite{BQ1} and adapted by Goyt and Sagan~\cite{SagGoytqFib} to produce some identities involving some of these new $q$-Fibonacci numbers.  In Sections 5, we consider the cycle decomposition of the permutations in $S_n(123,132,213)$ and determine two different $q$-Fibonacci numbers from these.  Finally, in Section 6, we study the distribution of {\it maj} and {\it inv} over pattern avoiding sets of permutations that West~\cite{West} showed are counted by the even Fibonacci numbers.

\section{Distributions and $q$-Fibonacci Numbers}

Let $s(\pi)$ be the number of singletons of $\pi$ and $d(\pi)$ be the number of doubletons of $\pi$.  We let $$F_n^I(x,y,q)=\sum_{\pi\in S_n(123,132,213)}x^{s(\pi)}y^{d(\pi)}q^{inv(\pi)},$$ and $$F_n^{I'}(x,y,q)=\sum_{\pi\in S_n(231,312,321)}x^{s(\pi)}y^{d(\pi)}q^{inv(\pi)}.$$  Also let $$F_n^M(x,y,q)=\sum_{\pi\in S_n(123,132,213)}x^{s(\pi)}y^{d(\pi)}q^{maj(\pi)},$$ and $$F_n^{M'}(x,y,q)=\sum_{\pi\in S_n(231,312,321)}x^{s(\pi)}y^{d(\pi)}q^{maj(\pi)}.$$

Let the {\it block structure} of a (reverse) layered permutation be a word in the set $A=\{s,d\}^*$, where the $k^{th}$ letter of the word is an $s$ ($d$) if the $k^{th}$ layer is a singleton (doubleton).  For example, the block structure of the permutation $\pi=6753421$ is the word $v_{\pi}=dsdss$.  It's not hard to see that a (reverse) layered permutation is uniquely defined by its block structure.

If $v$ is a word in $A$ then let the {\it length} of $v$, $\ell(v)$, be the sum of the lengths of its letters, where $\ell(s)=1$ and $\ell(d)=2$.  For example, $\ell(dsdss)=7$.  Let $A_n=\{v\in A:\ell(v)=n\}$.  There is an obvious bijection $\phi:S_n(123,132,213)\rightarrow A_n$.  For any letter $a$ of a word $v$ let $a_{lv}$ ($a_{rv}$) be the subword of $v$ consisting of the letters to the left (right) of $a$ in $v$.  

We define two weights on these words as follows.  Let $v=a_1a_2\dots a_n\in A$, and let $\omega_i(v)=\omega_i(a_1)\cdot\omega_i(a_2)\cdots\omega_i(a_n)$, where $\omega_i(s)=xq^{\ell(s_{rv})}$ and $\omega_i(d)=yq^{2\ell(d_{rv})}$.  Similarly, let $\omega_m(v)=\omega_m(a_1)\cdot\omega_m(a_2)\cdots\omega_m(a_n)$, where $\omega_m(s)=xq^{\ell(s_{lv})}$, and $\omega_m(d)=yq^{\ell(d_{lv})}$.  

Using the running example $\pi=6753421$ and $v_\pi=dsdss$, we have $$\omega_i(\phi(\pi))=x^3y^2q^{19}=x^{s(\pi)}y^{d(\pi)}q^{inv(\pi)}$$ and $$\omega_m(\phi(\pi))=x^3y^2q^{16}=x^{s(\pi)}y^{d(\pi)}q^{maj(\pi)}.$$ Thus, we may redefine our $q$-Fibonacci numbers in the following way, $$F_n^I(x,y,q)=\sum_{\pi\in S_n(123,132,213)}\omega_i(\phi(\pi)),$$ and $$F_n^M(x,y,q)=\sum_{\pi\in S_n(123,132,213)}\omega_m(\phi(\pi)).$$

\begin{thm}  $F_0^I(x,y,q)=1$, $F_1^I(x,y,q)=x$, and for $n\geq2$, $$F_n^I(x,y,q)=xq^{n-1}F_{n-1}^I(x,y,q)+yq^{2(n-2)}F_{n-2}^I(x,y,q).$$\end{thm}

{\bf Proof:}  $A_0$ consists of one word with no doubletons and no singletons, so $F_0^I(x,y,q)=1$.  $A_1$ consists of the word $s$, and the weight of this word is $x$.  Each word in $A_n$ begins with $s$ or $d$, whose weight is $\omega_i(s)=xq^{n-1}$ and $\omega_i(d)=yq^{2(n-2)}$ respectively.  All but the first letter in the word is a word in $A_{n-1}$ or $A_{n-2}$ respectively.  This gives us the identity.  \bbox

\medskip

\begin{thm}  $F_0^M(x,y,q)=1$, $F_1^M(x,y,q)=x$, and for $n\geq2$, $$F_n^M(x,y,q)=xq^{n-1}F^M_{n-1}(x,y,q)+yq^{n-2}F^M_{n-2}(x,y,q).$$\end{thm}

{\bf Proof:}  As above $A_0$ gives us $F_0^M(x,y,q)=1$, and $A_1$ gives us that $F_1^M(x,y,q)=x$.  Each word in $A_n$ ends with $s$ or $d$, whose weight is $\omega_m(s)=xq^{n-1}$ and $\omega_m(d)=yq^{n-2}$ respectively.  All but the last letter in the word is in $A_{n-1}$ or $A_{n-2}$ respectively.  This proves the identity.  \bbox

\medskip

The next two Lemmas explain how $F^I(x,y,q)$ is related to $F^{I'}(x,y,q)$ and how $F^M(x,y,q)$ is related to $F^{M'}(x,y,q)$.

\begin{lem}  For $n\geq0$,
$$F_n^I(x,y,q)=q^{n\choose2}F_n^{I'}\left(x,y,\frac{1}{q}\right).$$
\end{lem}

{\bf Proof:}  The left hand side is the distribution of $inv$ on $S_n(123,132,213)$.  On the right, $F_n^{I'}(q)$ is the distribution of $inv$ on $S_n(231,312,321)$.  Recall that $S_n(123,132,213)$ is the set of reverse layered matchings and $S_n(231,312,321)$ is the set of layered matchings.  Let $\pi\in S_n(123,132,213)$, and $$\rho:S_n(123,132,213)\rightarrow S_n(231,312,321)$$  be defined by $\rho(\pi)=\bar{\pi}$, the reversal of $\pi$.  Then $\pi$ and $\bar{\pi}$ have the same number of doubletons, say $k$.  It's not hard to see that $inv(\pi)={n\choose2}-k$, and $inv(\bar{\pi})=k.$  Thus, $$q^{inv(\pi)}=q^{n\choose2}\cdot\left(\frac{1}{q}\right)^{inv(\bar{\pi})}.$$
\bbox

\medskip

\begin{lem}  For $n\geq0$,
$$F_n^M(x,y,q)=q^{n\choose2}F_n^{M'}\left(x,y,\frac{1}{q}\right).$$
\end{lem}

{\bf Proof:}  The left hand side is the distribution of $maj$ on $S_n(123,132,213)$.  On the right, $F_n^{M'}(q)$ is the distribution of $maj$ on $S_n(231,312,321)$.  Let $\pi\in S_n(123,132,213)$, and $$\psi:S_n(123,132,213)\rightarrow S_n(231,312,321)$$  be defined by $\psi(\pi)=\tilde{\pi}$ where $\pi$ and $\tilde{\pi}$ have the same block structure.  Then descents in $\tilde{\pi}$ take place in the positions where descents do not take place in $\pi$, and vice versa.  Thus, $$q^{maj(\pi)}=q^{n\choose2}\cdot\left(\frac{1}{q}\right)^{maj(\tilde{\pi})}.$$
\bbox

\section{$F_n^M(x,y,q)$ and Previous $q$-Fibonacci Numbers}

The recursion found in Theorem 2.2 is the same recursion found by Goyt and Sagan~\cite{SagGoytqFib} for their $q$-Fibonacci number $F_n(x,y,q)$, which involves the $rb$ statistic.  We will give a bijection from $S_n(123,132,213)$ to $\Pi_n(13/2,123)$ (defined below) that maps the {\it maj} statistic to the $rb$ statistic.  In order to discuss this bijection, we must first talk about pattern avoidance in set partitions.

A partition $\alpha$ of $[n]$, denoted $\alpha\vdash[n]$, is a family of disjoint subsets $B_1,\:B_2,\dots,\:B_k$ of $[n]$, called {\it blocks}, such that $\bigcup_{i=1}^kB_i=[n]$, and $B_i\not=\emptyset$ for each $i$.  We write $\alpha=B_1/B_2/\dots /B_k$, omitting set braces and commas, and we always list the blocks in the standard order, where $$\min B_1<\min B_2<\dots<\min B_k,$$ and the elements in each block are in ascending order.

A {\it layered} partition is a partition of the form $\alpha=[1,i]/[i+1,j]/\dots/[k+1,n]$, and a {\it matching} is a partition $B_1/B_2\dots /B_k$, where $|B_i|\leq 2$ for $1\leq i\leq k$.  For example, $12/345/6/78$ is a layered partition, and $1/23/4/5/67/89$ is a layered matching.  As before, one element blocks will be called singletons and two element blocks will be called doubletons.

Suppose $\alpha=A_1/A_2/\dots /A_k\vdash[m]$ and $\beta=B_1/B_2/\dots/B_\ell\vdash[n]$.  We say that $\alpha$ is contained in $\beta$, $\alpha\subseteq\beta$, if there are distinct blocks $B_{i_1},\:B_{i_2},\dots,\:B_{i_k}$ of $\beta$, such that $A_j\subseteq B_{i_j}$.  For example, if $\beta=1/236/45$ then $\alpha'=26/4$ is contained in $\beta$, but $\alpha'=1/2/3$ is not because the 2 and the 3 would have to be in separate blocks of $\beta$.

Suppose $\alpha=A_1/A_2/\dots/A_k\vdash[m]$ and $\beta=B_1/B_2/\dots /B_\ell\vdash[n]$.  We say $\beta$ contains the {\it pattern} $\alpha$ if there is some $\alpha'\subseteq \beta$ such that $\alpha'$ and $\alpha$ are order isomorphic, otherwise we say that $\beta$ avoids $\alpha$.  

Define $$\Pi_n=\{\alpha\vdash[n]\}$$ and for any set of partitions $R$, $$\Pi_n(R)=\{\alpha\vdash[n]:\pi\mbox{ avoids every partition in $R$}\}.$$

Goyt~\cite{Goytpartitions3} showed that all permutations in the set $\Pi_n(13/2,123)$ are layered matchings and are counted by the Fibonacci numbers.  Like the layered permutations, each $\alpha\in\Pi_n(13/2,123)$ is uniquely determined by its block structure.

We now turn our attention to the $rb$ statistic developed by Wachs and White~\cite{WachsWhite}.  Let $\alpha=B_1/B_2/\dots/B_k$ be a partition and $b\in B_i$.  Then $(b,B_j)$ is a {\it right bigger pair} of $\alpha$ if $j>i$ and $\max B_j>b$.  For example, in the partition $\alpha=1/236/45$, $(3,\{4,5\})$ is a right bigger pair.  Let $rb(\alpha)$ be the number of right bigger pairs in $\alpha$.  We will say that the block $B_i$ contributes $t$ to the $rb$ statistic if there are $t$ right bigger pairs of the form $(b,B_i)$.  It is immediately apparent from the definition of layered matchings that the contribution of a block $B_i$ in a layered matching is $\min B_i-1$.  Let $$F_n(x,y,q)=\sum_{\alpha\in\Pi_n(13/2,123)}x^{s(\alpha)}y^{d(\alpha)}q^{rb(\alpha)}$$ be the $q$-Fibonacci number associated with the $rb$ statistic.  

We will now define a weight on the block structure of a layered matching.  Let $\alpha$ be a layered matching and let $v_\alpha=a_1a_2\dots a_k$ be its block structure, then we define $\omega_{rb}$ to be $\omega_{rb}(v)=\omega_{rb}(a_1)\cdot\omega_{rb}(a_2)\cdots\omega_{rb}(a_k)$, where $\omega_{rb}(s)=xq^{\ell(s_{lv})}$ and $\omega_{rb}(d)=yq^{\ell(d_{lv})}$.  For example, if $\alpha=12/3/4/56/78$, then its block structure is $v_\alpha=dssdd$ and $\omega_{rb}(v_\alpha)=x^2y^3q^{15}=x^{s(\alpha)}y^{d(\alpha)}q^{rb(\alpha)}$.  Let $\eta:S_n(123,132,213)\rightarrow\Pi_n(13/2,123)$, where $\pi$ and $\eta(\pi)$ have the same block structure.  By the definition of $\omega_m$ and $\omega_{rb}$ we have that $maj(\pi)=rb(\eta(\pi))$.

\begin{thm} For $n\geq0$, $$F_n^M(x,y,q)=F_n(x,y,q).$$\end{thm} \bbox

\medskip

In his paper~\cite{Ciglerqfib1} Cigler describes Morse sequences, which relate to our $q$-Fibonacci polynomials $F_n^{M'}(x,y,q)$.  A {\it Morse sequence} of length $n$ is a sequence of dots and dashes, where each dot has length $1$ and each dash has length $2$.  For example, $v=\bullet\bullet--\bullet-$ is a Morse sequence of length $9$.  Let $MS_n$ be the set of Morse sequences of length $n$.  Each Morse sequence corresponds to a layered matching where a dot is replaced by a singleton block and a dash by a doubleton.  So, $|MS_n|=F_n$.

Let $\mu=m_1m_2\dots m_k$ and let $\varphi:MS_n\rightarrow S_n(231,312,321)$ satisfy $\varphi(\mu)=\pi$, where $\pi$ has $k$ blocks and block $i$ is a singleton if $m_i$ is a dot or a doubleton if $m_i$ is a dash.  Clearly, $\varphi$ is a bijection.

Cigler defines the weight of a dot to be $0$ and the weight of a dash to be $a+1$ where $a$ is the length of the portion of the sequence appearing before the dash.  Also, he lets $w(\mu)$ be the sum of the weights of the dashes in $\mu$.  For example, the sequence above has weight $3+5+8=16$.  Let 
$$F_n^C(x,y,q)=\sum_{\mu\in MS_n}x^{t(\mu)}y^{h(\mu)}q^{w(\mu)},$$
where $t(\mu)$ is the number of dots and $h(\mu)$ is the number of dashes in $\mu$.  In~\cite{Ciglerqfib1}, Cigler shows that $F_n^C(x,y,q)$ satisfies $F_0^C(x,y,q)=1$, $F_1^C(x,y,q)=x$, and $$F_n^C(x,y,q)=xF_{n-1}^C(x,y,q)+yq^{n-1}F_{n-2}^C(x,y,q).$$

\begin{lem}

The map $\varphi$ described above satisfies for any $\mu\in MS_n$, $$w(\mu)=maj(\varphi(\mu)).$$

\end{lem}

{\bf Proof:}  Let $\mu=m_1m_2\dots m_k\in MS_n$.  Let $\pi=p_1p_2\dots p_n\in S_n(231,312,321)$ such that $\varphi(\mu)=\pi$.  Note that if $m_j=-$, then the corresponding block in $\pi$ is a doubleton.  In $S_n(231,312,321)$, descents only take place in the first position of the doubletons.  If $m_j=-$, and $m_j$ contributes $k$ to $w(\mu)$, then the length of the sequence before $-$ is $k-1$.  Thus, $m_j$ corresponds to the doubleton $p_kp_{k+1}$, which contributes $k$ to $maj(\pi)$.  If $m_j=\bullet$, then $m_j$ contributes $0$ to $w(\mu)$, and the corresponding singleton $p_k$ is not the beginning of a doubleton and contributes $0$ to $maj(\pi)$.  \bbox

\medskip

The following theorem is an immediate consequence of Lemma 3.2, so we omit its proof.

\begin{thm}  For $n\geq0$,
$$F_n^{M'}(x,y,q)=F_n^C(x,y,q).$$
\end{thm}
\bbox

\section{Inversion Theorems}

We now turn our attention to Fibonacci identities and give bijective proofs of identities involving $F_n^I(x,y,q)$.  These proofs use the same techniques of Benjamin and Quinn~\cite{BQ1}, and Goyt and Sagan~\cite{SagGoytqFib}.

\begin{thm} For $m,n\geq1$,
$$F_{m+n}^I(x,y,q)=F_m^I(xq^n,yq^{2n},q)F_n^I(x,y,q)+yq^{2(n-1)}F_{m-1}^I(xq^{n+1},yq^{2(n+1)},q)F_{n-1}^I(x,y,q).$$
\end{thm}

{\bf Proof:}  Let $\pi=p_1p_2\dots p_{m+n}\in S_{m+n}(123,132,213)$, and suppose $p_mp_{m+1}$ is not a doubleton. Then $v_\pi=v'v''$ where $v'$ is a word in $S_m(123,132,213)$ and $v''$ is a word in $S_n(123,132,213)$.  Since there are $n$ elements to the right of $v'$, the weight of each singleton in $v'$ is increased by a factor of $q^n$ and each doubleton by $q^{2n}$.  Thus, we get the first part of our identity.

Now suppose $p_mp_{m+1}$ is a doubleton.  In this case $v_\pi=v'dv''$ where $v'$ is a word in $S_{m-1}(123,132,213)$ and $v''$ is a word in $S_{n-1}(123,132,213)$.  Since there are $n+1$ elements to the right of $v'$, the weight of each singleton in $v'$ is increased by a factor of $q^{n+1}$ and each doubleton by $q^{2(n+1)}$.  The doubleton $p_mp_{m+1}$ has weight $yq^{2(n-1)}$.  Thus, we get the second part of the identity.

Clearly, all permutations fall into one of these two cases, so we obtain the desired identity.  \bbox

Setting $m=1$ in the previous identity leads to the identity in Theorem 2.1.  It's also interesting to note that if we set $n=1$ then we obtain the identity $F_m(x,y,q)=xF_{m-1}(xq,yq^2,q)+yF_{m-2}(xq^2,yq^4,q)$ for $m\geq 2$.  This identity may be obtained in the same way that the identity in Theorem 2.1 was obtained except that we divide $A_n$ into two sets by whether the words in $A_n$ end in a singleton or a doubleton.  

\begin{thm}  For $n\geq0$,
$$F_n^I(x,y,q)=\sum_{2k\leq n}{n-k\choose k}x^{n-2k}y^kq^{{n\choose 2}-k}.$$
\end{thm}

{\bf Proof:}  Let $\pi\in S_n(123,132,213)$.

There is exactly one permutation with no doubletons.  In this case, $v_\pi=ss\dots s$ and $\displaystyle \omega_i(v_\pi)=x^nq^{\sum_{j=0}^{n-1}j}=x^nq^{n\choose2}.$

Consider the set of words in $v_\pi$ with exactly $k$ $d$'s.  There are $n-k$ letters and therefore ${n-k\choose k}$ such words.  

Notice that for each $d$ in $v_\pi$, the power of $q$ in $ss\dots s$ is reduced by one.  Thus, every word, $v_\pi$, with exactly $k$ doubletons satisfies $\omega_i(v_\pi)=x^{n-2k}y^kq^{{n\choose 2}-k}$.  Summing over all possible $k$ gives us the desired identity.  \bbox

The classical Fibonacci polynomials are $F_n(x,y)={n-k\choose k}x^{n-2k}y^k=F_n^I(x,y,1)$.  There are many identities involving classical Fibonacci polynomials, and it turns out that we can translate most of these into identities involving $F_n^I(x,y,q)$.  To do this we will need two other identities.

\begin{thm}  For $n\geq0$, $$F_n^I(xq,yq^2,q)=q^nF_n^I(x,y,q),$$ and $$F_n^I(x,y,q)=q^{n\choose 2}F_n^I\left(x,\frac{y}{q},1\right).$$\end{thm}

{\bf Proof:}  For the first identity place a phantom 1 at the end of every word in $A_n$ and increase each element by 1.  This would increase the weight of every singleton by one and every doubleton by two.  On the other hand, the singleton would be involved in $n$ inversions.  

The proof of the second identity is essentially the same as the proof of the previous theorem.  Each doubleton reduces the maximum number of inversions, ${n\choose 2}$, by one.  \bbox

\medskip

The well known Cassini identity is $F_n(x,y)^2-F_{n+1}(x,y)F_{n-1}(x,y)=(-1)^ny^n$.  This can be translated into a Cassini-like identity for $F_n^I(x,y,q)$ using the second identity from Theorem 4.3.  The first thing we do is replace $y$ by $\frac{y}{q}$ in the identity above and obtain $F_n(x,\frac{y}{q},1)^2-F_{n+1}(x,\frac{y}{q},1)F_{n-1}(x,\frac{y}{q},1)=(-1)^n(\frac{y}{q})^n$.  Now, multiply through by $q^{(n^2-n+1)}$ and obtain $$\left(q^{n\choose2}F_n(x,\frac{y}{q},1)\right)^2-q^{n+1\choose 2}F_{n+1}(x,\frac{y}{q},1)q^{n-1\choose2}F_{n-1}(x,\frac{y}{q},1)=(-1)^ny^nq^{(n-1)^2}.$$  Using the second identity from Theorem 4.3, we obtain a the Cassini-like identity for $F_n^I(x,y,q)$ as follows $$qF_n(x,y,q))^2-F_{n+1}(x,y,q)F_{n-1}(x,y,q)=(-1)^ny^nq^{(n-1)^2}.$$

The following theorems give more bijective proofs of $q$-Fibonacci identities involving $F^I(x,y,q)$.  

\begin{thm} For $n\geq0$, $$F_{n+2}^I(x,y,q)=x^{n+2}q^{{n+2\choose2}}+\sum_{j=0}^nx^{n-j}yq^{\frac{n^2+3n-j^2+j}{2}}F_j^I(x,y,q).$$
\end{thm}

{\bf Proof:}   Let $\pi=p_1p_2\dots p_{n+2}\in S_{n+2}(123,132,213).$  There is exactly one such permutation with no doubletons.  Thus, $v_\pi=ss\dots s$ and $\omega_i(v_\pi)=x^{n+2}q^{\sum_{j=1}^{n+1}j=x^{n+2}q^{n+2\choose2}}$.

Let the first doubleton in $\pi$ be $p_{n-j+1}p_{n-j+2}$.  Thus, $v_\pi=ss\dots sdv'$, and we have that $\omega_i(v_\pi)=\omega_i(ss\dots s)w_i(d)w_i(v')$.  Notice that $\omega_i(ss\dots s)=x^{n-j}q^{\sum_{k=j+2}^{n+1}k}=x^{n-j}q^{{n+2\choose 2}-{j+2\choose 2}}$, and $\omega_i(d)=yq^{2j}$.  Thus, \begin{eqnarray*}\omega_i(sss\dots s)w_i(d)&=&x^{n-j}yq^{{n+2\choose 2}-{j+2\choose 2}+2j}\\ &=&x^{n-j}yq^{\frac{n^2+3n-j^2+j}{2}}.\end{eqnarray*}
We have that $\omega_i(v')$ contributes $F_j^I(x,y,q)$.  Summing over $j$ gives the desired identity.  \bbox

\begin{thm} For $n\geq0$,
$$F_{2n+1}^I(x,y,q)=\sum_{j=0}^nxy^jq^{4nj-2j^2+2n-2j}F_{2n-2j}^I(x,y,q).$$
\end{thm}

{\bf Proof:}  Let $\pi\in S_{2n+1}(123,132,213)$.  Since $2n+1$ is always odd, every permutation must contain at least one singleton.  Assume there are $j$ doubletons to the left of the first singleton.  Then $v_\pi=ddd\dots dsv'$.  Thus, $\omega_i(v_\pi)=\omega_i(ddd\dots d)\omega_i(s)\omega_i(v')$.    We can see that $\displaystyle \omega_i(ddd\dots d)=y^jq^{\sum_{k=1}^{j}{2(2n-2k+1)}}=y^jq^{4nj-2j^2}$.  Also, $\omega_i(s)=xq^{2n-2j}$, and $\omega_i(v')$ contributes $F_{2n-2j}^I(x,y,q)$.  Summing over $j$ gives the identity. \bbox

\begin{thm} For $n\geq0$,
$$F_{2n}^I(x,y,q)=y^nq^{n(n-1)}+\sum_{j=0}^{n-1}xy^jq^{4nj-2j^2-4j+2n-1}F_{2n-2j-1}^I(x,y,q).$$
\end{thm}

{\bf Proof:}  Let $\pi\in S_{2n}(123,132,213)$.

Since $2n$ is even, there is exactly one such permutation with no singletons.  So we have $v_\pi=ddd\dots d$ and $ \omega_i(v_\pi)=y^nq^{\sum_{k=0}^{n-1}2k}=y^nq^{n(n-1)}$, which gives us the first term.

Let the first singleton be $p_{2j+1}$.  Then $v_\pi=ddd\dots dsv'$ where $v'$ is a word in $A_{2n-(2j+1)}$. So $\omega_i(v_\pi)=\omega_i(ddd\dots d)\omega_i(s)\omega_i(v')$.  Then $\omega_i(ddd\dots d)= y^jq^{\sum_{k=1}^j2(2n-2k)}=y^jq^{4nj-2j^2-2j}$, and $\omega_i(s)=xq^{2n-(2j+1)}$.  Summing over $j$ gives the desired identity. \bbox

\begin{thm} For $n\geq0$,
$$F_{n+1}^I(x,y,q)F_n^I(x,y,q)=\sum_{j=0}^nxy^{n-j}q^{(n-j)(n+j-1)+j}\left(F_j^I(x,y,q)\right)^2.$$
\end{thm}

{\bf Proof:}  Let $(\pi_1,\pi_2)\in S_{n+1}(123,132,213)\times S_n(123,132,213)$, $v_{\pi_1}=a_1a_2\dots a_k$, and $v_{\pi_2}=b_1b_2\dots b_\ell$.  We search through the words in the order $a_1,b_1,a_2,b_2,\dots$ until we find the first $s$.  This will happen because either $n$ or $n+1$ is odd.  

Suppose the first $s$ is some $a_i$.  Then $v_{\pi_1}=ddd\dots dsv'$.  Assume there are $\frac{n-j}{2}$ doubletons to the left of $s$, where $n-j$ is even.  Then $\omega_i(v_{\pi_1})=\omega_i(ddd\dots d)\omega_i(s)\omega_i(v')$.  So $\omega_i(s)=xq^j$ and $\omega_i(v')$ contributes $F_j^I(x,y,q)$.  We can also see that $v_{\pi_2}=ddd\dots dv''$, and $v_{\pi_2}$ also begins with $\frac{n-j}{2}$ doubletons.  Thus $\omega_i(v_{\pi_2})=\omega_i(ddd\dots d)\omega_i(v'')$, with $\omega_i(v'')$ contributing $F_j^I(x,y,q)$.  Thus, the weight of the doubletons at the beginning of $\pi_1$ and $\pi_2$ is $\displaystyle y^{n-j}q^{\sum_{k=1}^{n-j}2(n-k)},$ which is $y^{n-j}q^{(n-j)(n+j-1)}$.  Thus, $\omega_i(v_{\pi_1})\omega_i(v_{\pi_2})$ contributes $xy^{n-j}q^{(n-j)(n+j-1)}\left(F_j^I(x,y,q)\right)^2$.

Suppose the first $s$ is some $b_i$.  Then $v_{\pi_2}=ddd\dots dsv'$.  Assume there are $\frac{n-j-1}{2}$ doubletons to the left of $s$, where $n-j$ is odd.  Then $\omega_i(v_{\pi_2})=\omega_i(ddd\dots d)\omega_i(s)\omega_i(v'')$.  So $\omega_i(s)=xq^j$ and $\omega_i(v'')$ contributes $F_j^I(x,y,q)$.  We can also see that $v_{\pi_1}=ddd\dots dv'$, and that $v_{\pi_1}$ begins with $\frac{n-j+1}{2}$ doubletons.  Thus $\omega_i(v_{\pi_1})=\omega_i(ddd\dots d)\omega_i(v')$, with $\omega_i(v')$ contributing $F_j^I(x,y,q)$.  So the weight of the doubletons at the beginning of $\pi_1$ and $\pi_2$ is again $y^{n-j}q^{\sum_{k=1}^{n-j}2(n-k)}=y^{n-j}q^{(n-j)(n+j-1)}.$  Again, we must have that $\omega_i(v_{\pi_1})\omega_i(v_{\pi_2})$ contributes $xy^{n-j}q^{(n-j)(n+j-1)}\left(F_j^I(x,y,q)\right)^2$.

Summing over all possible $j$ gives the desired identity.

\section{Cycle Decomposition}

We now turn our attention to two cycle decomposition statistics.  Recall that a permutation can be decomposed into cycles.  For example, the permutation $\sigma=978645312$ has cycle decomposition $(192738)(465)$.  We are interested in the distribution of two different statistics.  Let $c(\sigma)$ be the number of cycles of $\sigma$ and $c_i(\sigma)$ to be the number of cycles of $\sigma$ of length $i$, where the length of a cycle is the number of elements in the cycle.  Define
$$F_n^D(x,y,q)=\sum_{\sigma\in S_n(123,132,213)}x^{s(\sigma)}y^{d(\sigma)}q^{c(\sigma)}.$$
and
$$F_n^{D'}(x,y,z_1,z_2,z_3,\dots)=\sum_{\sigma\in S_n(123,132,213)}x^{s(\sigma)}y^{d(\sigma)}\prod_{i\geq1}z_i^{c_i(\sigma)}.$$

To determine recursive forms for these two Fibonacci polynomials we will need to be able to determine the cycle decomposition from the word associated to $\sigma\in S_n(123,132,213)$.  

For example, the permutation $\sigma=978645312$, with cycle decomposition $(192738)(465)$, gives us the word $sdsdsd$.  The first cycle in the decomposition arises from $\underline{9}\underline{78}645\underline{3}\underline{12}$, corresponding to the letters  $\underline{s}\underline{d}sd\underline{s}\underline{d}$.  Notice that this cycle is produced by jumping back and forth between the beginning and end of the permutation.  This suggests that we should do the same when using the word to produce a cycle.  

Consider the word $v_\sigma=sdsdsd$ from above.  We will rewrite this word by alternating taking a letter from the beginning of the word and the end of the word.  For example, the new word associated to $v_\sigma$ is $v_\sigma'=sddssd$ because the first letter in $v_\sigma$ is an $s$, the last (sixth) letter is a $d$, the second letter is a $d$, the fifth letter is an $s$, etc.  The first cycle of $\sigma=(192738)(465)$ comes from the first four letters of $v_\sigma'$, namely $sdds$, and the second cycle comes from the last two letters of $v_\sigma'$, namely $sd$.

Given a permutation $\sigma$, its associated word $v_\sigma$ and its new word $v_\sigma'$, we can use $v_\sigma'$ to determine the lengths of the cycles and number of cycles in an inductive way.  Each $v_\sigma'$ begins with a word that will always produce a certain type of cycle.  These {\it prefixes} are $ss,dd,ds*s,dssd,sd^\ell s,$ and $dsd^\ell s$, where the $*$ in $ds*s$ may represent $s$ or $d$.  We'll explain what happens to the cycle decomposition when our word begins with one of these prefixes.

If $v_\sigma'$ begins with the prefix $ss$ then $\sigma$ is of the form $n\dots 1$.  It is easy to see that the resulting cycle is $(1n)$, which is of length $2$.  If $v_\sigma'$ begins with $dd$ then $\sigma$ is of the form $(n-1)n\dots12$.  Thus, we have $(1(n-1))(2n)$, which is two cycles of length $2$.  The prefix $ds*s$ gives us that $\sigma$ is of the form $(n-1)n\dots21$.  Thus, we have $(1(n-1)2n)$, which is a cycle of length 4.    The prefix $dssd$ gives us that $\sigma$ is of the form $(n-1)n(n-2)\dots231$.  Thus, we have $(1(n-1)3(n-2)2n)$, which is a cycle of length 6.  

The next theorem and corollary take care of the other two cases.

\begin{thm} The prefix $sd^\ell s$ gives a $2\ell+2$ cycle.
\end{thm}

{\bf Proof:}  Suppose without loss of generality that $\ell$ is even.  If the encoding of a permutation $\sigma$ looks like $sd^\ell s$ then 
$$\sigma=n(n-2)(n-1)\dots(n-\ell)(n-(\ell-1))\dots(\ell+1)(\ell-1)\ell\dots12.$$
By simply reading through these elements of the permutation we get the cycle
$$(1n2(n-2)4\dots\ell(n-\ell)(\ell+1)(n-(\ell-1))(\ell-1)\dots(n-1)).$$
The case when $\ell$ is odd is proved similarly. \bbox

\medskip

\begin{cor}  The prefix $dsd^\ell s$ gives a $2\ell+4$ cycle.
\end{cor}
\bbox

\medskip

It is easy to see that all possible prefixes are accounted for.  We may now use these prefixes to produce a couple of $q$-Fibonacci identities based on cycle decomposition statistics.

\begin{thm}  For $n\geq0$, $$F_{n+2}^D(x,y,q)=x^2qF_n^D(x,y,q)+\left(y^2q^2+2x^2yq\right)F_{n-2}^D(x,y,q)+2\sum_{k=3}^{\left\lfloor{\frac{n}{2}}\right\rfloor}x^2y^{k-1}qF_{n-2k}^D(x,y,q).$$
\end{thm}

{\bf Proof:}  If the prefix has length 2, then the first blocks are $ss$, which contributes $x^2q$.  The remaining blocks are counted by $F_n^D(x,y,q)$.  When the prefix is of length $4$, the block structure is $dd$, $sds$, or $ds*s$.  For $dd$, the first cycle contributes $y^2q$.  For $sds$ and $ds*s$, the first cycle contributes $x^2yq$.  The remaining blocks are counted by $F_{n-2}^D(x,y,q)$.  When the prefix is of length $6$, the block structure is $dssd$ or $sdds$.  Both contribute $x^2y^2q$, with the remaining blocks counted by $F_{n-6}^D(x,y,q)$.

When the prefix is length $2\ell$, with $4\leq \ell\leq\left\lfloor\frac{n}{2}\right\rfloor$, the block structure is $sd^{\ell-1}s$ or $dsd^{\ell-2}s$.  The first cycle then contributes $x^2y^{\ell-1}q$.  The remaining blocks are counted by $F_{n-2\ell}^D(x,y,q)$.  The desired identity is achieved. \bbox

\medskip

The proof of the following theorem is similar to the proof above, so we omit it.

\begin{thm} For $n\geq0$, $$F_{n+2}^{D'}(x,y,z_1,z_2,z_3,\dots)=x^2z_2F_n^{D'}(x,y,z_1,z_2,z_3,\dots)+\left(y^2z_2^2+2x^2yz_4\right)F_{n-2}^{D'}(x,y,z_1,z_2,z_3,\dots)$$$$+2\sum_{k=3}^{\left\lfloor{\frac{n}{2}}\right\rfloor}
x^2y^{k-1}z_{2k}F_{n-2k}^{D'}(x,y,z_1,z_2,z_3,\dots).$$
\end{thm}
\bbox

\section{Even $q$-Fibonacci Numbers}

In~\cite{West}, West uses generating trees to show that certain sets of permutations avoiding a pattern from $S_3$ and a pattern from $S_4$ are counted by even Fibonacci numbers.  In this section we will describe a few $q$-Fibonacci numbers that arise from studying the distribution of $inv$ on some of these sets.  

The first set that we will consider is the set $S_n(123,2143)$.  Let $$F_{2n-2}^{W_1}(q)=\sum_{\sigma\in S_n(123,2143)}q^{inv(\sigma)}.$$  We will start with a basic identity that arises directly from the construction described by West.

\begin{thm} We have $F_1^{W_1}(q)=1$ and for $n\geq 2$, $$F_{2n}^{W_1}(q)=q^{n-1}F_{2n-2}^{W_1}(q)+\sum_{k=2}^{n}q^{(n-1)(k-1)+{k\choose2}}F_{2(n-k)}^{W_1}(q).$$ \end{thm}

{\bf Proof:}  We will use West's description of how permutations in $S_n(123,2143)$ are constructed in order to produce this identity.  We will construct a permutation in $S_n(123,2143)$ from a permutation $\sigma=p_1p_2\dots p_{n-1}\in S_{n-1}(123,2143)$ by placing $n$ immediately before $p_k$, $1\leq k\leq n$, and avoiding copies of 123 and 2143.  If $n$ is placed before $p_k$ then we say $n$ is placed in the $k^{th}$ gap.  If $k=n$ then $n$ is placed at the end of the permutation.  

We can place $n$ in the first or second gap without producing a copy of 123 or 2143.  Suppose $\sigma$ begins with $p_1\leq n-2$.  If we want to place $n$ in the $k^{th}$ gap where $k\geq 3$, we must have $p_1p_2\dots p_{k-1}$ descending to avoid 123.  However, since $p_1\leq n-2$, we must have $n-1$ appearing after $n$, and $p_1p_2n(n-1)$ gives a copy of 2143.  This means that if $p_1\leq n-2$, then $n$ may only be placed in the first or second gap.

The first term of the identity is given by placing $n$ in the first gap.  If $n$ is placed in the $k^{th}$ gap $k\geq 3$, then $p_1p_2\dots p_{k-1}=(n-1)(n-2)\dots(n-k+1)$, and $p_k\dots p_{n-1}$ must be a permutation in $S_{n-k}(123,2143)$.  The contribution of $p_1p_2\dots p_{k-1}n$ is $q^{\sum_{\ell=1}^{k-1} n-1-\ell}=q^{(n-1)(k-1)+{k\choose2}}$.  If we sum over all possible permutations $p_k\dots p_n$ in $S_n(123,2143)$ we obtain $F_{2(n-k)}(q)$.  Summing over $2\leq k\leq n$ and including the term given by placing $n$ in the first gap gives the desired identity.  \bbox

\medskip

The next set we will consider is the set $S_n(132,3241)$.  Let 
$$F_{2n-2}^{W_2}(q)=\sum_{\sigma\in S_n(132,3241)}q^{inv(\sigma)}.$$

\begin{thm}  $F_0^{W_2}(q)=1$ and for $n\geq1$, 
$$F_{2n-2}^{W_2}(q)=(q^{n-1}+1)F_{2n-4}^{W_2}(q)+\sum_{k=1}^{n-2}q^{k(n-k)}F_{2n-2k-4}^{W_2}(q).$$
\end{thm}

{\bf Proof:}  We construct a permutation in $S_n(132,3241)$ from $\sigma\in S_{n-1}(132,3241)$ by placing $n$ in one of the $k$ gaps of $\sigma$.  It's not hard to see that one can place $n$ in the first or $n^{th}$ gap without creating a copy of $132$ or $3241$.  This gives the first term in the identity.

Suppose now that $n$ is not in the first or $n^{th}$ gap.  We can see that all elements to the right of $n$ must be smaller than all elements to the left of $n$ in order to avoid $132$.  Also, all elements to the left of $n$ must be in ascending order to avoid $3241$.  If $n$ is in the $k^{th}$ gap, then $\sigma=(n-k+1)(n-k+2)\dots np_kp_{k+1}\dots p_{n-1}$.   The first $k+1$ elements of $\sigma$ contributes $q^{k(n-k)}$ and if we sum over all possible permutations $p_k\dots p_{n-1}$ gives $F_{2n-2k-4}(q)$.  Summing over $2\leq k\leq n-1$ and including the term given by placing $n$ in the first or $n^{th}$ gap gives the desired identity. \bbox

\medskip

The next set we will consider is the set $S_n(132,3412)$.  Let 
$$F_{2n-2}^{W_3}(q)=\sum_{\sigma\in S_n(132,3412)}q^{inv(\sigma)}.$$

\begin{thm}  $F_0^{W_3}(q)=1$ and for $n\geq1$, 
$$F_{2n-2}^{W_3}(q)=(q^{n-1}+1)F_{2n-4}^{W_2}(q)+\sum_{k=1}^{n-2}q^{k(n-k)+{n-k\choose2}}F_{2k-4}^{W_3}(q).$$
\end{thm}

{\bf Proof:}  The proof of this theorem is very similar to the above, only in this case if we put $n$ in the $k^{th}$ gap where $k\not=1$ or $n$, then all the elements to the right must be smaller than all of the elements to the left, and the elements to the right must be descending.   \bbox

\medskip

\bibliography{Bib1}

\end{document}